\documentclass[12]{amsart}

\newtheorem{theorem}{Theorem}[section]
\newtheorem{lemma}[theorem]{Lemma}

\theoremstyle{definition}

\theoremstyle{remark}

\numberwithin{equation}{section}
\begin{document}
\title{On minimal norms on $M_n$}
\author[M. Mirzavaziri, M. S. Moslehian]{Madjid Mirzavaziri$^1$ and Mohammad Sal
Moslehian$^2$}
\address{Department of Mathematics, Ferdowsi University, P. O. Box 1159, Mashhad 91775, Iran;
\newline Banach Mathematical Research Group (BMRG), Mashhad, Iran.}
\email{mirzavaziri@math.um.ac.ir and madjid@mirzavaziri.com}
\address{$^2$ Department of Mathematics, Ferdowsi University, P. O. Box 1159, Mashhad 91775, Iran;
\newline Centre of Excellence in Analysis on
Algebraic Structures (CEAAS), Ferdowsi University, Iran.}
\email{moslehian@ferdowsi.um.ac.ir and moslehian@ams.org}
\subjclass[2000]{Primary 15A60; Secondary 47A30, 46B99}
\keywords{Induced norm; generalized induced norm; algebra norm; the
full matrix algebra; minimal norm.}

\begin{abstract}
In this note, we show that for each minimal norm $N(\cdot)$ on the
algebra $M_n$ of all $n \times n$ complex matrices, there exist
norms $\|\cdot\|_1$ and $\|\cdot\|_2$ on ${\mathbb C}^n$ such that
$$N(A)=\max\{\|Ax\|_2: \|x\|_1=1,~ x\in {\mathbb C}^n\}$$
for all $A \in M_n$. This may be regarded as an extension of a known
result on characterization of minimal algebra norms.
\end{abstract}
\maketitle


\section{Introduction}

Let ${\mathcal M}_n$ denote the algebra of all $n\times n$ complex
matrices $A$ with entries in ${\mathbb C}$ together with the usual
matrix operations. By an algebra norm (or a matrix norm) we mean a
norm $\|\cdot\|$ on ${\mathcal M}_n$ such that
$\|AB\|\leq\|A\|\|B\|$ for all $A, B\in {\mathcal M}_n$. It is easy
to see that the norm
$\|A\|_{\sigma}=\displaystyle{\sum_{i,j=1}^n}|a_{ij}|$ is an algebra
norm, but the norm $\|A\|_{m}=\max\{|a_{i,j}|: 1\leq i,j\leq n\}$ is
not an algebra norm, see \cite{H-M-M}.

\noindent Let $\|\cdot\|_1$ and $\|\cdot\|_2$ be two norms on
${\mathbb C}^n$. Then the norm $\|\cdot\|_{1,2}$ on ${\mathcal
M}_n$ defined by $\|A\|_{1,2}:=\max\{\|Ax\|_2: \|x\|_1=1\}$ is
called the generalized induced (or g-ind) norm constructed via
$\|\cdot\|_{1}$ and $\| . \|_{2}$. If $\|\cdot\|_1=\|\cdot\|_2$,
then $\|\cdot\|_{1,1}$ is called an induced norm.

\noindent It is known that
$\|A\|_C=\max\{\displaystyle{\sum_{i=1}^n}|a_{i,j}|: 1\leq j\leq
n\}, \|A\|_R=\max\{\displaystyle{\sum_{j=1}^n}|a_{i,j}|: 1\leq
i\leq n\}$ and the spectral norm $\|A\|_S= \max\{\sqrt{\lambda}:
\lambda \textrm{~is~ an ~eigenvalue~ of~} A^*A\}$ are induced by
$\ell_1, \ell_\infty$ and $\ell_2$, respectively; cf. \cite{H-J}.
Recall that the $\ell_p$-norm ($1\leq p\leq \infty$) on ${\mathbb
C}^n$ is defined by
$$\ell_p(x)=\ell_p(\displaystyle{\sum_{i=1}^n}x_ie_i)=\left
\{ \begin{array}{cc}(\displaystyle{\sum_{i=1}^n}|x_i|^p)^{1/p}&1\leq p<\infty\\
\max\{|x_1|, \ldots, |x_n|\}&p=\infty \end{array}\right .$$

\noindent It is known that the algebra norm $\|A\|=\max\{\|A\|_C,
\|A\|_R\}$ is not induced and it is not hard to show that it is
not g-ind too; cf. Corollary 3.2.6 of \cite{B-L}.

\noindent A norm $N(\cdot)$ on $M_n$ is called minimal if for any
norm $|||\cdot|||$ on $M_n$ satisfying $|||\cdot||| \leq N(\cdot)$
we have $|||\cdot|||= N(\cdot)$. It is known \cite[Theorem
3.2.3]{B-L} that an algebra norm is an induced norm if and only if
it is a minimal element in the set of all algebra norms. Note
that a generalized induced norm may not be minimal. For instance,
put $\|\cdot\|_\alpha=\ell_{\infty}(.),
\|\cdot\|_\beta=2\ell_{2}(.)$ and $\|\cdot\|_\gamma=\ell_{2}(.)$.
Then $\|\cdot\|_{\gamma,\beta}\leq\|\cdot\|_{\alpha,\beta}$ but
$\|\cdot\|_{\gamma,\beta}\neq\|\cdot\|_{\alpha,\beta}$.

\noindent In \cite{H-M-M}, the authors investigate generalized
induced norms. In particular, they examine the problem that ``for
any norm $\|\cdot\|$ on $M_n$, there are two norms $\|.\|_1$ and
$\|.\|_2$ on ${\mathbb C}^n$ such that $\|A\|=\max\{\|Ax\|_{2}:
\|x\|_{1}=1\}$ for all $A\in {\mathcal M}_n$?'' In this short note,
we utilize some ideas of \cite{H-M-M} to study the minimal norms on
${\mathcal M}_n$. More precisely, we show that for each minimal norm
$N(\cdot)$ on the algebra $M_n$ of all $n \times n$ complex
matrices, there exists norms $\|\cdot\|_1$ and $\|\cdot\|_2$ on
${\mathbb C}^n$ such that $N(A)=\max\{\|Ax\|_2: \|x\|_1=1,~ x\in
{\mathbb C}^n\}$ for all $A \in M_n$. In particular, if $N(\cdot)$
is an algebra norm, then $\|.\|_1=\|\cdot\|_2$. This may be regarded
as an extension of the above known result on characterization of
minimal algebra norms.


\section{Main Result}

For $x \in {\mathbb C}^n$ and $1\leq j\leq n$, let $C_{x,j}\in
{\mathcal M}_n$ be defined by the operator $C_{x,j}(y)=y_{j}x$.
Hence $C_{x,j}$ is the $n\times n$ matrix with $x$ in the $j$ column
and $0$ elsewhere. Define $C_{x}\in {\mathcal M}_n$ by
$C_{x}=\sum_{j=1}^nC_{x,j}$. Hence $C_{x}$ is the $n\times n$ matrix
whose all columns are $x$.

\noindent If $\|\cdot\|_{1,2}$ is a generalized induced norm on
${\mathcal M}_n$ obtained via $\|\cdot\|_1$ and $\|\cdot\|_2$ then
$\|C_x\|_{1,2}=\alpha\|x\|_2$, where $\alpha=\max\{|\sum_{j=1}^n
y_j|: \|(y_1,\ldots,y_j,\ldots,y_n)\|_1=1\}$.


\noindent To achieve our goal, we need the following lemmas.

\begin{lemma}\label{lem1} \cite[Theorem 2.7]{H-M-M}
Let $\|\cdot\|_1$ and $\|\cdot\|_2$ be two norms on ${\mathbb C}^n$.
Then $\|\cdot\|_{1,2}$ is an algebra norm on ${\mathcal M}_n$ if and
only if $\|\cdot\|_1\leq\|\cdot\|_2$.
\end{lemma}


\begin{lemma}\label{lem2} \cite[Corollary 2.5]{H-M-M}
$\|\cdot\|_{1,2}=\|\cdot\|_{3,4}$ if and only if there exists
$\gamma>0$ such that $\|\cdot\|_1=\gamma\|\cdot\|_3$ and
$\|\cdot\|_2=\gamma\|\cdot\|_4$.
\end{lemma}


\begin{theorem}
Let $N(\cdot)$ be a minimal norm on ${\mathcal M}_n$, then
$N(\cdot)=\|\cdot\|_{1,2}$ for some $\|.\|_1$ and $\|.\|_2$ on
${\mathbb C}^n$. Moreover, if $N(\cdot)$ is an algebra norm, then
$\|.\|_1=\|\cdot\|_2$.
\end{theorem}
\begin{proof}
For $x\in {\mathbb C}^n$, set $$\|x\|_1=\max\{N(C_{Ax}): N(A)=1,
A\in {\mathcal M}_n\}$$ and
$$\|x\|_2=N(C_x).$$

\noindent We shall show that $\|\cdot\|_1$ and $\|\cdot\|_2$ are
norms on ${\mathbb C}^n$.

\noindent To see $\|\cdot\|_1$ is a norm, let $x \in {\mathbb C}^n$.
Then $\|x\|_1=0$ if and only if $N(C_{Ax})=0$ for all matrix $A$
with $N(A)=1$, and this holds if and only if $Ax=0$ for all $A$, or
equivalently $x=0$.

\noindent For $\alpha \in {\mathbb C}^n$ and $x, y \in {\mathbb
C}^n$, we  have
\begin{eqnarray*}
\|\alpha x\|_1 &=& \max\{N(C_{A(\alpha x)}): N(A)=1, A\in {\mathcal
M}_n\}\\
&=& \max\{N(\alpha C_{Ax}): N(A)=1, A\in {\mathcal
M}_n\}\\
&=& \max\{|\alpha| N(C_{Ax}): N(A)=1, A\in {\mathcal
M}_n\}\\
&=& |\alpha| \max\{N(C_{Ax}): N(A)=1, A\in {\mathcal
M}_n\}\\
&=& |\alpha|~\|x\|_1\,.
\end{eqnarray*}
and
\begin{eqnarray*}
\|x+y\|_1 &=& \max\{N(C_{A(x+y)}): N(A)=1, A\in {\mathcal
M}_n\}\\
&=& \max\{N(C_{Ax}+C_{Ay}): N(A)=1, A\in {\mathcal
M}_n\}\\
&\leq& \max\{N(C_{Ax}): N(A)=1, A\in {\mathcal M}_n\}\\
&& + \max\{N(C_{Ay}): N(A)=1, A\in {\mathcal
M}_n\}\\
&=& \|x\|_1 + \|y\|_1\,.
\end{eqnarray*}
\noindent To see $\|\cdot\|_2$ is a norm, let $x \in {\mathbb C}^n$.
Then $\|x\|_2=0$ if and only if $C_x=0$ and this holds if and only
if $x=0$.

\noindent For $\alpha \in {\mathbb C}^n$ and $x, y \in {\mathbb
C}^n$, we have
\begin{eqnarray*}
\|\alpha x\|_2 = N(C_{\alpha x}) = N(\alpha C_x)= |\alpha| N(C_x)
= |\alpha|~\|x\|_2\,.
\end{eqnarray*}
and
\begin{eqnarray*}
\|x+y\|_2 = N(C_{x+y})= N(C_x+C_y) \leq N(C_x) + N(C_y) = \|x\|_2 +
\|y\|_2\,.
\end{eqnarray*}

Now let $A\in {\mathcal M}_n\backslash\{0\}$. Then
$N(\frac{A}{N(A)})=1$. So that
$$\|\frac{A}{N(A)}(x)\|_2 =\|N(C_{\frac{A}{N(A)}(x)}) \leq \|x\|_1$$
whence
$$\|Ax\|_2 \leq N(A)\|x\|_1\,.$$
Therefore $\|A\|_{1,2} \leq N(A)$. Since $N(\cdot)$ is a minimal
norm, we conclude that $\|A\|_{1,2} = N(A)$.

\noindent If $N(A)$ is an algebra norm, then Lemma \ref{lem1}
implies that $\|\cdot\|_1 \leq \|\cdot\|_2$.

\noindent Next let $A \in {\mathcal M}_n$. It follows from $\|Ax\|_1
\leq \|A\|_{1,1}\|x\|_1 \leq \|A\|_{1,1}\|x\|_2, \quad (x \in
{\mathbb C}^n)$ that $\|A\|_{2,1} \leq \|A\|_{1,1}$. In a similar
fashion one can get
\begin{eqnarray*}
\|\cdot\|_{2,1} \leq \|\cdot\|_{k,k} \leq \|\cdot\|_{1,2}\qquad
(k=1,2).
\end{eqnarray*}
By the minimality of $\|\cdot\|_{1,2}$ we deduce that
$\|\cdot\|_{1,2}=\|\cdot\|_{1,1}$. It then follows from Lemma
\ref{lem2} that $\|\cdot\|_1=\|\cdot\|_2$.
\end{proof}



\begin{thebibliography}{99}
\bibitem{B-L} G. R. Belitski\u\i and Yu. I. Lyubich, \textit{Matrix Norms and their Applications},
Translated from the Russian by A. Iacob.: Operator Theory: Advances
and Applications, 36. Birkhäuser Verlag, Basel, 1988.
\bibitem{H-M-M} S. Hejazian, M. Mirzavaziri and M. S. Moslehian, \textit{Generalized induced norms}, Czechoslovak Math. J. \textbf{57} (2007), no. 1, 127--133.
\bibitem{H-J} R. A. Horn and C. R. Johnson, \textit{Matrix Analysis},
Cambridge University Press, Cambridge, 1994.
\end{thebibliography}
\end{document}